\documentclass[sn-mathphys-num]{sn-jnl}

\usepackage{graphicx}%
\usepackage{multirow}%
\usepackage{amsmath,amssymb,amsfonts}%
\usepackage{amsthm}%
\usepackage{mathrsfs}%
\usepackage[title]{appendix}%
\usepackage{xcolor}%
\usepackage{textcomp}%
\usepackage{manyfoot}%
\usepackage{booktabs}%

\usepackage{xspace}                                             
\usepackage{amstext}                                            
\usepackage{array}                                              
\usepackage{bbm}                                                
\usepackage{subcaption}                                         
\usepackage{tikz,bbding}                                        
\usetikzlibrary{shapes.geometric, positioning, arrows, calc}    
\usetikzlibrary {arrows.meta}                                   

\usepackage[textwidth=30mm,
  backgroundcolor=yellow,
  linecolor=orange,
  textsize=tiny]{todonotes}
\makeatletter 
\@mparswitchfalse%
\makeatother
\normalmarginpar 

\newcolumntype{C}{>{$}c<{$}}                                    
\newcolumntype{L}{>{$}l<{$}}                                    
\newcolumntype{R}{>{$}r<{$}}                                    
\DeclareMathOperator*{\argmin}{argmin}                          
\def\1{{\mathbbm{1}}}                                           
\def\R{{\mathbb{R}}}                                            
\def\RU{\overline{\R}}                                          
\newcommand{\solar}{{\sf solar}\xspace}                         
\newcommand{\nomad}{{\sf NOMAD}\xspace}                         
\newcommand{\mads}{{MADS}\xspace}                               
\newcounter{algocounter}                                        
\renewcommand{\thealgocounter}{3.\arabic{algocounter}}          
\newenvironment{algo}                                           
    {\refstepcounter{algocounter}}
    {}
\newtheorem{assumption}{Assumption}                             

\newtheorem{theorem}{Theorem}
\newtheorem{lemma}[theorem]{Lemma}

\begin{document}

\title[Fidelity and interruption control for expensive constrained multi-fidelity blackbox optimization]{Fidelity and interruption control for expensive constrained multi-fidelity blackbox optimization}


\author[2,3]{\fnm{St\'ephane} \sur{Alarie}}\email{alarie.stephane@hydroquebec.com (0000-0001-6693-4509)}

\author[1,2]{\fnm{Charles} \sur{Audet}}\email{charles.audet@gerad.ca (0000-0002-3043-5393)}

\author[3]{\fnm{Miguel} \sur{Diago}}\email{diagomartinez.miguel3@hydroquebec.com (0000-0001-7857-8070)}

\author[1,2]{\fnm{S\'ebastien} \sur{Le~Digabel}}\email{sebastien.le.digabel@gerad.ca (0000-0003-3148-5090)}

\author*[1,2]{\fnm{Xavier} \sur{Lebeuf}}\email{xavier.lebeuf@polymtl.ca (0000-0001-6945-7349)}

\affil[1]{\orgdiv{Department of Mathematical and Industrial Engineering}, \orgname{Polytechnique Montr\'eal}, \orgaddress{\street{2500 Chem. de Polytechnique}, \city{Montr\'eal}, \postcode{H3T~1J4}, \state{Qu\'ebec}, \country{Canada}}}

\affil[2]{\orgname{GERAD}, \orgaddress{\street{2920 Chem. de la Tour}, \city{Montr\'eal}, \postcode{H3T~1N8}, \state{Qu\'ebec}, \country{Canada}}}

\affil[3]{\orgname{Hydro-Qu\'ebec}, \orgaddress{\street{1800 Bd Lionel-Boulet}, \city{Varennes}, \postcode{J3X~1S1}, \state{Qu\'ebec}, \country{Canada}}}


\abstract{
This work introduces a novel blackbox optimization algorithm for computationally expensive constrained multi-fidelity problems.
When applying a direct search method to such problems, the scarcity of feasible points may lead to numerous costly evaluations spent on infeasible points.
Our proposed fidelity and interruption controlled optimization algorithm addresses this issue by leveraging multi-fidelity information, allowing for premature interruption of an evaluation when a point is estimated to be infeasible.
These estimations are controlled by a biadjacency matrix, for which we propose a construction.
The proposed method acts as an intermediary component bridging any non multi-fidelity direct search solver and a multi-fidelity blackbox problem, giving the user freedom of choice for the solver.
A series of computational tests are conducted to validate the approach. 
The results show a significant improvement in solution quality when an initial feasible starting point is provided.
When this condition is not met, the outcomes are contingent upon specific properties of the blackbox.

}

\keywords{
Blackbox optimization,
Derivative-free optimization,
Multi-fidelity,
Constrained optimization,
Direct search methods,
Static surrogates
}

\maketitle

\section{Introduction}
\label{sec:intro}

This work considers the constrained optimization problem
\begin{equation*}
    \qquad \qquad\min_{x\in \Omega} \ f(x)\quad\text{ where }\quad 
    \Omega=\{x\in X:c_j(x)\leq 0, j\in J\}
    \tag{$\mathcal{P}$}
    \label{eq-pb}
\end{equation*}
in which~$X = [\ell, u] \subset \R^n$ is a set defined by
    unrelaxable constraints with $\ell, u \in \R^n$,
    $f:X \rightarrow\RU = \R \cup \{\infty\}$
    and 
    $c_j:X\rightarrow\RU$, $j \in J=\{1,2,\dots,m\}$,
    are the objective and quantifiable constraint functions, respectively. 
The set of feasible points~$\Omega$, is delimited by the constraint functions~$c_j(x)\leq 0, j\in J$ and by the bounds $[\ell, u]$.
These constraints form the vector~$c(x)=(c_1(x),c_2(x),\dots,c_m(x))$.
We use~$\RU$ because~$f$ can be set to~$\infty$ to reject points when using a barrier method,
    and constraints can be set to~$\infty$ when a hidden constraint is violated.

The present work considers~$f$ and~$c$ as functions for which derivatives are unavailable,
    and as outputs of a multi-fidelity blackbox~\cite{AuHa2017},
    which is expensive to run and may fail to evaluate.
Fidelity is defined as the degree to which a model reproduces the state and behavior of the true object, 
    represented here by the real scalar~$\phi$, an element of a finite discrete subset of~$[0,1]$.
An evaluation of~$f$ and~$c$ at~$\phi=1$ results in the highest precision,
    and usually the highest cost.
Conversely,
    an evaluation at~$\phi<1$ may be interpreted as the evaluation of a static surrogate model.
    
The computational time required to evaluate the blackbox at a trial point~$x\in X$ using the fidelity~$\phi$ is denoted by~$t(x,\phi) \in \R_+$. The functions~$f$,~$c$ and~$c_j$ appearing in Problem~\eqref{eq-pb} are expanded to~$f(x,\phi)$,~$c(x,\phi)$, and~$c_j(x,\phi)$, with~$f(x)=f(x,1)$,~$c(x)=c(x,1)$, and~$c_j(x)=c_j(x,1)$ for $j\in J$, where the parameter~$\phi$ indicates the fidelity of an evaluation. It is not assumed that with~$\phi_a<\phi_b$,~$c(x,\phi_a)$ and~$f(x,\phi_a)$ can be deduced from a blackbox evaluation at fidelity~$\phi_b$.

Recent literature concerning multi-fidelity blackbox optimization problems predominantly emphasizes research in the unconstrained case. 
The present work focuses on the exploitation of constraint information from multi-fidelity.
This research serves as a first step, to propose an algorithmic approach that comprehensively leverages the impact of fidelity on both constraint and objective function values.

\subsection{Motivation}


The study of multi-fidelity is motivated
    by an asset management blackbox optimization problem
    encountered at Hydro-Qu\'ebec as part of the PRIAD project~\cite{PRIAD_CoBLALDeKoMe2020,GaChKoCoHeBlDeAb2021,PRIAD_KoMeCoGaVoAlDeBl2021},
    which is constituted of  computationally expensive blackboxes wherein the violation of some constraints can be predicted with low fidelity evaluations.
Since the overall time allowed to the optimization process is limited,
    strategies need to be devised to accelerate the evaluations.
    
One of the first strategy that comes to mind is parallelism, as 
    in~\cite{Alba13,GrKo06,Haftka16},
    but is not sufficient to solve the problem in the allowed time.
As discussed in~\cite{PRIAD_KoMeCoGaVoAlDeBl2021},
    it would require thousands of processors to solve the problem 
    within a month or less,
    even if one assumes linear improvement with the number of processors.
In addition,
    this assumption is unlikely to be satisfied since the speed gain diminishes when more processors are used in the computation~\cite{HiMa2008}.
Instead,
    we suggest to solve Problem~\eqref{eq-pb} based on
    (i) the preemption concept of~\cite{RaToMaThMaSe2010},
        which allows to interrupt an evaluation as soon as it is shown that
        a trial point~$x\in X$ will not replace the current incumbent solution,
    and 
    (ii) the idea of converging to a local solution by iteratively increasing
    the fidelity~\cite{PoWe06a,WePo2005}.
The purpose of using these mechanisms is to reduce the time spent
    on evaluating uninteresting points,
    thereby increasing the total number of evaluations
    and exploring the solution domain more intensively.
In other words,
    the present work attempts to only engage the minimal computational effort
    to reach better solutions
    within a predetermined time budget.

\subsection{Organization}

The proposed approach embeds interruption opportunities into an iterative fidelity evaluation process to reduce computational time.
For a given blackbox problem, a biadjacency matrix is constructed to select the fidelity level required to estimate the feasibility of each constraint.
This construction is based on a sub-problem, for which a fast solving method is presented.
Interruptions occur when a constraint is estimated to be violated at a point from low fidelity information.
Then, the point is assumed infeasible without incurring time-consuming higher fidelity evaluations.

This document is structured as follows.
Section~\ref{sec:rev_lit} presents a literature review of constraints management
    in a multi-fidelity blackbox optimization framework.
Section~\ref{sec:metho} presents an algorithm
    to solve Problem~\eqref{eq-pb}.
Details will then be provided to construct the biadjacency matrix by solving the related optimization sub-problem, 
    and to set the resulting interruption mechanism in order to reduce computational time.
Section~\ref{sec:results} shows how the algorithm performs on different benchmarking blackboxes involving a solar thermal power plant simulator.
Finally,
    Section~\ref{sec:discu} discusses the results.

\section{Literature Survey}
\label{sec:rev_lit}

This work addresses single objective blackbox optimization problems~\cite{AuHa2017}.
It also categorizes constraints based on the taxonomy of~\cite{LedWild2015}.
A constraint violation function~$h:\R^m\rightarrow\RU$, which stems from filter methods~\cite{AuDe04a,FlLe02a,FlLeTo06}, is introduced in~\cite{AuDe09a}:
\begin{equation*}
    h(x):=
    \begin{cases}
        \sum\limits_{j=1}^m(\max\{c_j(x),0\})^2 & \text{ if }x\in X \\
        \infty & \text{ otherwise.}
    \end{cases}
\end{equation*}

This function measures the level of violation of relaxable constraints.
When~$x\in\Omega$, the function satisfies $h(x)=0$, and~$h(x)>0$ otherwise.
A first method to deal with constraints is called the extreme barrier (EB),
    which is divided into two unconstrained minimization phases~\cite{AuHa2017}. 
The first is the feasibility phase, where~$\min_{x\in \R^n}h(x)$ is solved while disregarding the value of~$f$, until a feasible point is found. 
Then, the optimality phase takes place, where~$\min_{x\in \R^n}f_\Omega(x)$ is solved, in which $f_\Omega(x) = f_(x)$ when $x \in \Omega$ and $f_\Omega(x)$ is set to~$\infty$ otherwise. 
A more sophisticated approach for dealing with quantifiable constraints is the progressive barrier (PB)~\cite{AuDe09a}. 
It introduces a threshold~$h_{\max}^k\in\RU$, initialized at~$\infty$, that progresses towards~$0$ as the iteration counter $k$ increases. 
Any trial point~$x$ such that~$h(x)>h_{\max}^k$ is rejected from consideration.
Two incumbent points are updated at the end of each iteration~$k$: the feasible solution~$x$ with the lowest value of~$f(x)$, named~$x^{\text{feas}}$, and the infeasible solution, named~$x^{\text{inf}}$, with the lowest value of~$f(x)$ among the trial points satisfying~$h(x)\leq h^k_{\max}$. 
The PB explores around both incumbent solutions, and as~$h_{\max}^k$ decreases,~$x^{\text{inf}}$ approaches the feasible region. 
Because it is frequent that~$f(x^{\text{inf}})<f(x^{\text{feas}})$, this may lead to good feasible solutions. 
The progressive-to-extreme barrier (PEB)~\cite{AuDeLe10} combines the EB and PB approaches. Each constraint is initially treated by the PB, and as soon as it is satisfied by the incumbent solution, it becomes treated by the EB for the remainder of the optimization process.

To reduce the computational burden of an expensive blackbox, the {\em two-phase interruptible EB} algorithm~\cite{G-2021-65} is a version of the EB adapted for problems where the constraint values $c_j$ are sequentially computed through independent blackboxes.
This strategy exploits the fact that~$h(x)$ is the sum of non-negative terms.
During the evaluation of trial point~$x$, the constraint violation function~$h(x)$ value is calculated cumulatively. 
An evaluation is interrupted as soon as~$h(x)$ exceeds the constraint violation function value of the incumbent point. 
A second approach to the same problem is the {\em hierarchical satisfiability with EB} algorithm~\cite{G-2021-65}. 

It consists of solving a sequence of $m$ optimization problems. 
The objective function of problem $j$, for $j \in J$, 
    is to minimize $c_j(x)$ and is subject to
    $c_i(x) \leq 0$ for $i \in \{1,2, \ldots, j-1\}$.
Each optimization is stopped as soon a feasible point with a nonpositive objective function value is found.
The starting point of each problem but the first is the final solution of the previous one. 
A feasible solution in found when all problems are successfully solved.
These approaches replicate the idea of interrupting a simulation during its course when it becomes known that it will not contribute to the optimization process~\cite{RaToMaThMaSe2010}.

To the authors' knowledge,
    the term multi-fidelity often refers to the use of only two fidelities (one high and one low), or a few more.
Moreover, in the multi-fidelity setting, low-fidelity sources are used to guide further sampling of the high-fidelity source, either by finding promising regions, or by training a model in the context of Bayesian optimization or other uses of Gaussian processes, and in the unconstrained case. Alternatively, constraints can be considered with via a penalty in the objective function~\cite{AgRaKoBu2023}.
Reviews of the usage of such methods in the last few decades are provided in~\cite{Fe2023,PeWiGu2018}.  Additionally, Bayesian optimization can be expanded to utilize multi-task Gaussian processes~\cite{SeSnAd2013,LeBa2019} and solve multi-objective problems~\cite{QuJiQiLeAn2024,IrKaDo2024}.
In this work, multiple fidelity levels are instead exploited to reduce the overall cost of sampling the true blackbox, and in the context of direct search methods for constrained problems.
 This topic is the subject of very limited research. Herein, a novel algorithm that aims to fill this gap is introduced.
Moreover, the proposed approach selects relevant fidelity levels instead of assuming that a low-fidelity source is necessarily helpful.

In~\cite{SeKaSh2018}, a heuristic method is proposed for unconstrained problems where the objective function can be queried at a continuous range of fidelities~$[0,1]$. 
The global optimization with surrogate approximation of constraints (GOSAC) algorithm~\cite{MuWo2017} uses radial basis function surrogate models for constraints to solve problems where constraints are given by an expensive blackbox, but the objective function is easy to evaluate. In~\cite{AuCaJa2020a}, machine learning is used to guide a direct search algorithm when hidden, binary and unrelaxable constraints are present. Similarly,~\cite{MeMuSi2023} suggests a machine learning approach to predict the violation of hidden constraints, but it is not integrated in an optimization algorithm.

The fidelity of a model only indicates its predictive capability, but indicates nothing on the computational cost of the model. Both concepts are combined to describe adequacy of a model in~\cite{BaKo19}. 
They propose a framework to evaluate the model adequacy, and show its use with the \mads algorithm~\cite{AuDe2006}. 
The search step is used to select points that minimize the error induced by low cost models within a trust region. 
This framework expands the use of \mads to multi-disciplinary design optimization~\cite{BaKo20} and time-dependent multidisciplinary design optimization~\cite{BaKo20a}.

The computational results in this paper are conducted with version~4 of the \nomad software package~\cite{nomad4paper}, which is an implementation of the \mads algorithm~\cite{AuDe2006}. When the multi-fidelity aspect of a blackbox is of stochastic nature, many variations on the \mads algorithm and new direct search algorithms have been proposed to take into account stochastic noise~\cite{AlAuBoLed2019,G-2019-30,AudIhaLedTrib2016,dzahini22}.

In the multi-fidelity literature, it is often the case that benchmarks used to assess the performance of new techniques are analytical in nature; that is, it is assumed the outputs come from a blackbox, but instead they come from a known mathematical expression~\cite{andres2022bifidelity}. In the industrial setting however, it is possible for the amount of low-fidelity sources to be virtually infinite, particularly when the output is acquired via a simulation which can be sped up, and therefore approximated. Hence, benchmarks such as~\cite{MScMLG,wang2017generic} should be favored.

\section{Exploiting multi-fidelity in expensive constrained problems}
\label{sec:metho}

This section describes the {\em fidelity and interruption controlled optimization algorithm}.
The method is composed of three steps to solve Problem~\eqref{eq-pb}.
First, a model of the constraints' feasibility when varying the fidelity is built through Latin hypercube (LH) sampling~\cite{McCoBe79a}.
Second, an optimization sub-problem which constructs a biadjacency matrix~$B$ is solved.
These two steps set up the interruption mechanism.
Finally, a direct search optimization process is launched using the interruption mechanism to reduce computational time. 

Before presenting the algorithm, 
    we first outline how the interruption mechanism is used in Section~\ref{sec:step3},
    and then explain how to build it in Section~\ref{sec:steps1et2}.
Next, Section~\ref{sec:reductions} provides details on solving the sub-problem.
Finally, the complete {\em fidelity and interruption controlled optimization algorithm} is presented in Section~\ref{sec:full_algo}.

\subsection{Direct search using multi-fidelity based interruptions}
\label{sec:step3}

An ordered discrete set of fidelities
\begin{equation*}
    \Phi = \{\phi_i \in [0,1]: i \in \{1,2,\ldots, L\}, 0 \leq \phi_1 < \phi_2 < \ldots < \phi_L = 1\}
\end{equation*}
with~$\phi_L = 1$ as its largest element is considered. 
In situations where the fidelities~$\phi$ are free to be chosen in the continuous interval~$[0,1]$, the set~$\Phi$ can be chosen to contain many elements, as unused fidelities will be automatically filtered out by the algorithm.
It may also only contain a few fidelities, each corresponding to a static surrogate model available to the user.

The fidelity value~$\phi=0$ indicates that only a priori constraints are evaluated.
Such constraints have a known explicit formulation, and do not require the execution of an expensive process~\cite{LedWild2015}. Inside a blackbox, a priori constraints are checked first and outputs are directly returned if any are violated.

The proposed method assigns each blackbox constraint to a fidelity in $\Phi$.
These assignments can be seen as vertices in a bipartite graph where the nodes of one part represent the~$L$ fidelities, and the nodes of the other represent the~$m$ constraints.
This information is encapsulated in a biadjacency matrix~$B\in\{0,1\}^{L\times m}$, where~$B_{ij}=1$ if constraint~$c_j\leq0$ is assigned to fidelity~$\phi_i$, and~$B_{ij}=0$ otherwise.
The assignment $B_{ij}=1$ indicates that $\phi_i$ is the smallest fidelity that may be trusted to predict if the constraint $c_j\leq0$ is satisfied or not.

Section~\ref{sec:steps1et2} describes how these assignments are computed.

\begin{figure}[ht!]
    \begin{center}
        \includegraphics[width=0.7\textwidth]{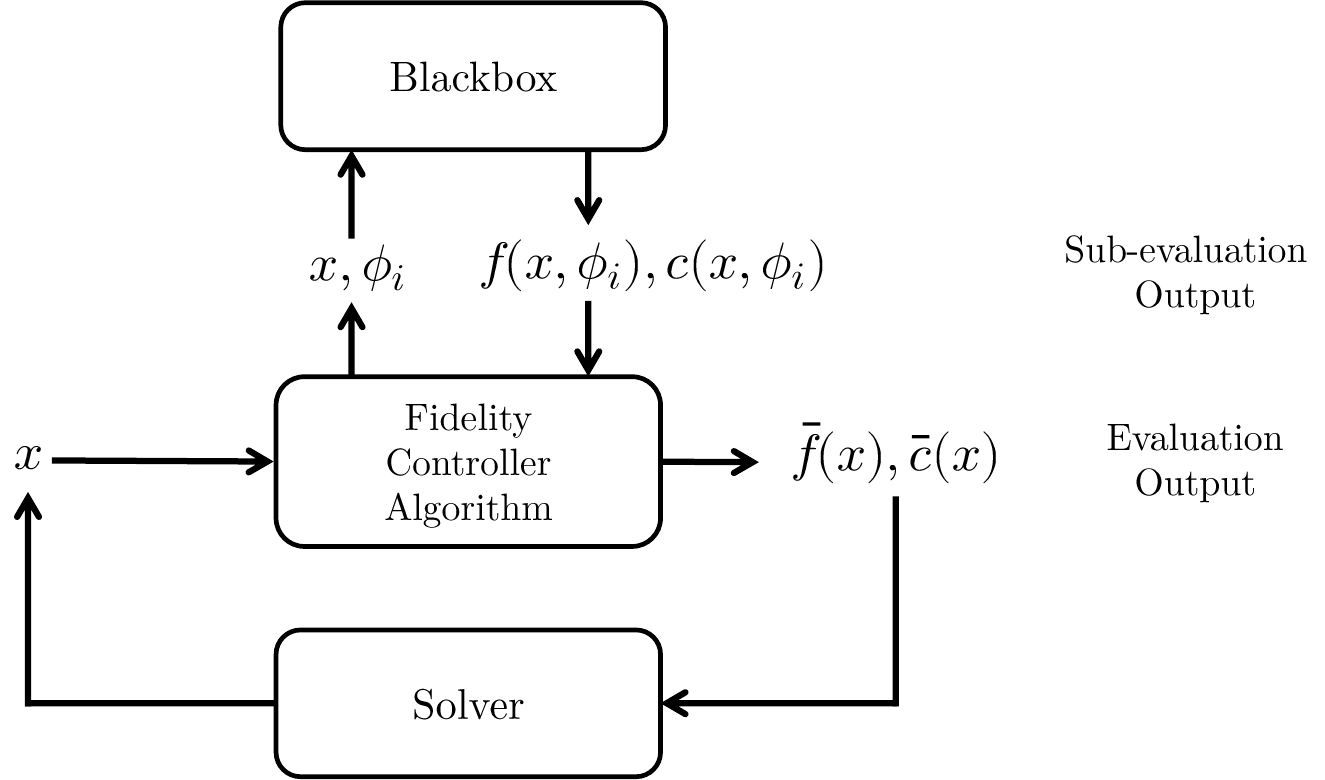}
    \end{center}
    \caption{Optimization loop with the {\em fidelity controller algorithm}. The notations~$\bar{f}(x)$ and~$\bar{c}(x)$ indicate the output values at~$x\in X$ at the last fidelity used by the algorithm.}
    \label{fig:optimization_loop}
\end{figure}

Multi-fidelity evaluations are used not to guide the exploration of the solution space, but rather to reduce overall computational time.
As such, the algorithm is coupled with a direct search solver that handles the exploration. When this solver determines that~$x$ is the next point to evaluate, the {\em fidelity controller algorithm} performs a sequence of calls, named sub-evaluations, to the blackbox in increasing order of fidelity in~$\Phi$ at~$x$, skipping the fidelities without assigned constraints.
The entire blackbox evaluation process at $x$ is interrupted as soon as a constraint~$c_j(x)\leq0$, for some $j \in J$, is violated at some fidelity~$\phi_i \in \Phi$, provided that the constraint is assigned to a fidelity of at most $\phi_i$.
Then, the most recent sub-evaluation's outputs are returned to the solver.
In this case, the solver may not receive the true blackbox output values, but it will still correctly deem~$x$ infeasible.
Figure~\ref{fig:optimization_loop} illustrates this process.

\begin{figure}[ht!]
    \centering
    \begin{tikzpicture}
        \node at (0.9,0.9) {Base case};
        
        \draw[draw] (0.2,0) rectangle ++(3,0.6);
        \draw[draw] (3.4,0) rectangle ++(3,0.6);
        \draw[draw] (6.6,0) rectangle ++(3,0.6);
        \draw[draw] (9.8,0) rectangle ++(3,0.6);
    
        \node[align=flush center,text width=0.3cm] at (0.2+1/2,0.3)     {$x^1$};
        \node[align=flush center,text width=0.3cm] at (3.4+1/2,0.3)     {$x^2$};
        \node[align=flush center,text width=0.3cm] at (6.6+1/2,0.3)     {$x^3$};
        \node[align=flush center,text width=0.3cm] at (9.8+1/2,0.3)     {$x^4$};
    
        \draw[draw] (0.2,-0.6) rectangle ++(3,0.6);
        \draw[draw] (3.4,-0.6) rectangle ++(3,0.6);
        \draw[draw] (6.6,-0.6) rectangle ++(3,0.6);
        \draw[draw] (9.8,-0.6) rectangle ++(3,0.6);
    
        \node[align=flush center,text width=0.3cm] at (0.2+1/2,-0.3)     {$\phi_3$};
        \node[align=flush center,text width=0.3cm] at (3.4+1/2,-0.3)     {$\phi_3$};
        \node[align=flush center,text width=0.3cm] at (6.6+1/2,-0.3)     {$\phi_3$};
        \node[align=flush center,text width=0.3cm] at (9.8+1/2,-0.3)     {$\phi_3$};
    \end{tikzpicture}
    \vspace{0.4cm}
    \begin{tikzpicture}
        \node at (1.5,0.9) {Fidelity controller};
        
    	\draw[draw] (0.2,0) rectangle ++(4.6,0.6);
        \draw[draw] (5  ,0) rectangle ++(0.5,0.6);
        \draw[draw] (5.7,0) rectangle ++(0.5,0.6);
        \draw[draw] (6.4,0) rectangle ++(1.6,0.6);
        \draw[draw] (8.2,0) rectangle ++(4.6,0.6);
    
        \node[align=flush center,text width=0.3cm] at (0.2+1/2,0.3)     {$x^1$};
        \node[align=flush center,text width=0.3cm] at (5  +1/4,0.3)     {$x^2$};
        \node[align=flush center,text width=0.3cm] at (5.7+1/4,0.3)     {$x^3$};
        \node[align=flush center,text width=0.3cm] at (6.4+1/2,0.3)     {$x^4$};
        \node[align=flush center,text width=0.3cm] at (8.2+1/2,0.3)     {$x^5$};
        
        \draw[draw] (0.2,-0.6) rectangle ++(0.5,0.6);
        \draw[draw] (0.7,-0.6) rectangle ++(1.1,0.6);
        \draw[draw] (1.8,-0.6) rectangle ++(3  ,0.6);
    
        \node[align=flush center,text width=0.3cm] at (0.2+1/4,-0.3)     {$\phi_1$};
        \node[align=flush center,text width=0.3cm] at (0.7+1/2,-0.3)     {$\phi_2$};
        \node[align=flush center,text width=0.3cm] at (1.8+1/2,-0.3)     {$\phi_3$};
    
        \draw[draw] (5  ,-0.6) rectangle ++(0.5,0.6);
        \draw[draw] (5.7,-0.6) rectangle ++(0.5,0.6);
        \draw[draw] (6.4,-0.6) rectangle ++(0.5,0.6);
        \draw[draw] (6.9,-0.6) rectangle ++(1.1,0.6);
    
        \node[align=flush center,text width=0.3cm] at (5  +1/4,-0.3)     {$\phi_1$};
        \node[align=flush center,text width=0.3cm] at (5.7+1/4,-0.3)     {$\phi_1$};
        \node[align=flush center,text width=0.3cm] at (6.4+1/4,-0.3)     {$\phi_1$};
        \node[align=flush center,text width=0.3cm] at (6.9+1/2,-0.3)     {$\phi_2$};
    
        \draw[draw] (8.2,-0.6) rectangle ++(0.5,0.6);
        \draw[draw] (8.7,-0.6) rectangle ++(1.1,0.6);
        \draw[draw] (9.8,-0.6) rectangle ++(3  ,0.6);
    
        \node[align=flush center,text width=0.3cm] at (8.2+1/4,-0.3)     {$\phi_1$};
        \node[align=flush center,text width=0.3cm] at (8.7+1/2,-0.3)     {$\phi_2$};
        \node[align=flush center,text width=0.3cm] at (9.8+1/2,-0.3)     {$\phi_3$};
    \end{tikzpicture}
    \begin{tikzpicture}[>=latex, scale=1, vertex/.style={draw, circle}]
        \node     (t) at (0   ,0)            {Time};
        \node     (0) at (0.35,0)            {};
        \draw[-{Latex[length=4mm]}] (0) -- ($(0)+(12,0)$);
    \end{tikzpicture}

    \caption{ An example showing that performing interruptions may save time compared to the base case that systematically uses the maximal fidelity.}
    \label{fig:fid_control_example}
\end{figure}
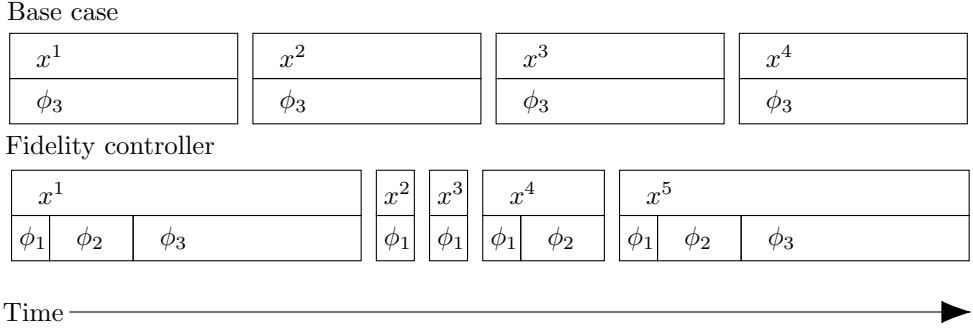

 Figure~\ref{fig:fid_control_example} shows an example with~$L=3$ that compares the {\em fidelity controller algorithm} with the base case that always uses the highest fidelity. 
The points~$x^k$ are generated by the solver, and are evaluated sequentially. 
The width of the rectangles represents evaluation time at a specified fidelity.
The widths are constant in the base case.
In this example, only~$x^1$ and~$x^5$ are feasible, making them longer to evaluate using the fidelity controller. 
However, the three other points are rapidly deemed infeasible and discarded.
This allows the solver to explore an additional point within the same time window.

Feasible points, which are evaluated at all fidelities in $\Phi$, take longer to compute with the {\em fidelity controller algorithm}.
In counterpart, prioritizing time-saving by assigning constraints to low fidelities carries the risk of misidentifying the feasibility of some points.
This implies that defining the biadjacency matrix requires careful consideration.

If it is allowed that no constraint is assigned to~$\phi=1$, meaning the true output values are never assessed, it is possible that the solution returned by the algorithm at the end of an optimization is said to be feasible but is actually not.
As this is highly undesirable, the method also considers~$f^*$, the value of the best feasible solution so far. 
Whenever a point~$x$ that has not been sub-evaluated with~$\phi=1$ is about to be labeled as the best solution so far by the solver (if~$f(x)<f^*$), an additional sub-evaluation with~$\phi=1$ is introduced to ensure the feasibility of~$x$.

Algorithm~\ref{algo:sous_eval_inter} shows the {\em fidelity controller algorithm}.
When values are returned before the algorithm's last line, the evaluation is said to be interrupted.
An optimization from a solver using Algorithm~\ref{algo:sous_eval_inter} constitutes the third step of the {\em fidelity and interruption controlled optimization algorithm}.

\begin{algo}
\label{algo:sous_eval_inter}
\begin{center}
    \begin{table}[ht]
    \centering \renewcommand{\arraystretch}{1.2}
    \begin{tabular}{L}
        \hline
        \textbf{Algorithm \thealgocounter: }\text{Fidelity controller algorithm} \\
        \hline
        \textbf{Input: } \text{ trial point } x; \text{ biadjacency matrix } B; \text{ incumbent value } f^*; \text{set of fidelities } \Phi \\
        \textbf{For each fidelity }\phi_i\in\Phi \text{ to which at least one constraint is assigned} \\
        \left.\middle|
                    \begin{array}{l}
                        \textbf{Evaluate } f(x,\phi_i) \text{ and } c(x,\phi_i)\text{, }\textbf{Store }\text{the values in $\bar f$ and $\bar c$.} \\
                        \textbf{If } c_j(x,\phi_i) > 0 \text{ for some } j\in J
                            \text{ where } B_{aj} = 1 \textbf{ for some } a \leq i \\
                        \left.\middle|
                                   \begin{array}{l}
                                        \textbf{Return }\bar f,\bar c
                                   \end{array}
                        \right. \\
                    \end{array}
        \right. \\
        \textbf{If }\text{no constraint is assigned to $\phi_L=1$ and }\bar f<f^* \\
        \left.\middle|
            \begin{array}{l}
                \textbf{Evaluate } f(x,1) \text{ and } c(x,1)\text{, }\textbf{Store }\text{the values in $\bar f$ and $\bar c$.}
            \end{array}
        \right. \\
    \textbf{Return }\bar f,\bar c \\
    \hline
    \end{tabular}
    \end{table}
\end{center}
\end{algo}

\subsection{Computing method for the biadjacency matrix}
\label{sec:steps1et2}

This section introduces a computing method for the matrix~$B$ which is broken down into two parts.
First, a model of how the constraints' feasibility is affected by the fidelity is proposed.
This process requires parallel computing to evaluate LH samples.
Second, this model is used to construct a~$B$ matrix.

\subsubsection{Constraint feasibility relative to fidelity model}

The user selects a sample size~$n_H$.
If a starting point~$x^0$ is provided, 
    the LH bounds are centered around~$x^0$ in order to model the constraints' feasibility where the optimization is likely to take place.
A new optimization parameter~$\rho\in[0,1]$ named the Latin hypercube sizing factor is introduced to indicate the size of the sampled region.
When~$\rho=1$, the region is equal to~$X = [\ell, u]$, and smaller values correspond to smaller regions.
The centered LH bounds~$(\ell^\text{cen},u^\text{cen})$ contain the following elements:
\begin{flalign}
    \ell^\text{cen}_i= & \,\text{max}\left(\ell_i,\ x^0_i-{\rho}(u_i-\ell_i)\right)\quad\forall\,i\in\{1,2,\dots,n\} \label{eq:licen}\\
    u^\text{cen}_i= & \,\text{min}\left(u_i,\ x^0_i+{\rho}(u_i-\ell_i)\right)\quad\forall\,i\in\{1,2,\dots,n\}. \label{eq:uicen}
\end{flalign}

If no starting point is provided, the LH bounds are those of~$X$. 
Let~$\text{H}\subset X$ denote the set of LH points,
    and let~$\Omega^{ap}$ be the set of points which do not violate any a priori constraint, where~$\Omega\subseteq\Omega^{ap}\subseteq X$.
Define~$I=\{1,2,\dots,L\}$ and the indicator function
\begin{equation*}
    \1(c_j(x,\phi)) :=\left\{
    \begin{array}{l}
        0\qquad\text{ if }c_j(x,\phi)\leq0 \\
        1\qquad\text{ otherwise}
    \end{array}
    \right. \\
    \qquad\forall j\in J.
\end{equation*}

A fidelity~$\phi_i \in \Phi$ is said to be {\em representative} for a constraint~$c_j\leq0$ at a point~$x$ if
\begin{equation*}
    \1(c_j(x,\phi))=\1(c_j(x,1))\quad\forall\phi\in\Phi\text{ with }\phi\geq\phi_i.
\end{equation*}
This definition allows to identify fidelities at which a sub-evaluation correctly identifies whether a constraint is violated or not on the LH samples. 
The behaviour of the constraints relative to fidelity is modeled in two ways: the likeliness of a fidelity to be representative for a given constraint at any point during the optimization, and the likeliness of a constraint to be violated at a given fidelity at any point during the optimization.
Another point of interest is the expected sub-evaluation time for the fidelities in~$\Phi$.
For any $i \in I$ and $j \in J$, define
\begin{alignat*}{2}
    r_{ij} & :=Pr[\text{fidelity $\phi_i$ is representative for constraint $c_j\leq0$ at 
 some }x\in X]  \\
    p_{ij} & :=Pr[c_j(x,\phi_i)>0 \mbox{ for some } x\in X]   \\
    t_i & :=\mathbb{E}_{x\in X}[t(x,\phi_i)].
\end{alignat*}
 Accordingly, the purpose of the LH samples is to estimate these statistical values. To achieve this, every sample is evaluated at each fidelity in~$\Phi$. For any $i \in I$ and $j \in J$:
\begin{alignat}{2}
    r_{ij} & \simeq \hat{r}_{ij}:=\frac{1}{\lvert\text{H}\cap\Omega^{ap}\rvert}\lvert\{x\in\text{H}\cap\Omega^{ap}:\phi_i\text{ is representative for }c_j\leq0\}\lvert\quad\label{eq:hatr}\\
    p_{ij} & \simeq \hat{p}_{ij}:=\frac{1}{\lvert\text{H}\cap\Omega^{ap}\rvert}\lvert\{x\in\text{H}\cap\Omega^{ap}:c_j(x,\phi_i)>0\}\lvert\label{eq:hatp}\\
    t_i & \simeq \hat{t}_i:=\frac{1}{\lvert\text{H}\cap\Omega^{ap}\rvert}\sum_{x\in\text{H}\cap\Omega^{ap}}t(x,\phi_i).\label{eq:hatt}
\end{alignat}
Note that when a point violates an a priori constraint, the other constraint values are unavailable (and uninteresting), hence the use of~$\text{H}\cap\Omega^{ap}$. This modeling of the constraints' behaviour relative to fidelity constitutes the first step of the {\em fidelity and interruption controlled optimization algorithm}.

\subsubsection{Optimal biadjacency matrix model}

A biadjacency matrix is computed based on given values of~$r_{ij}$,~$p_{ij}$ and~$t_i$. 
An explicit optimization Problem~$\mathcal{Q}$ is proposed below to compute~$B$, which minimizes the expected evaluation time during the optimization. It is also desired to minimize the probability of causing an interruption on a feasible point, but a bi-objective model is avoided by introducing a new threshold parameter~$\varepsilon\in[0,1]$, the upper bound on the probability that a constraint's feasibility is misidentified during the optimization.
Every element~$B_{ij}$ of the biadjacency matrix~$B$ corresponds to a decision variable in the model.
The additional decision variables~$y_i$ are defined as follows. They are continuous in~$\mathcal{Q}$ per their unimodality.
\begin{equation}
    y_i=
    \begin{cases}
        1 & \text{ if a sub-evaluation at fidelity~$\phi_i$ would be executed by Algorithm~\ref{algo:sous_eval_inter}} \\
        0 & \text{ otherwise.}
    \end{cases}
    \label{eq:def_y}
\end{equation}

In~$\mathcal{Q}$, the objective function~\eqref{eq:modele_obj} represents the expected evaluation time of a single point according to Algorithm~\ref{algo:sous_eval_inter}.
It is a sum of all the sub-evaluation times~$t_i$, multiplied by their probability of happening, which is the probability that no interruption happens earlier multiplied by~$y_i$ for each~$i\in I$.
It can be written as follows.
\begin{flalign}
    & t_1y_1 \nonumber\\
    & +t_2y_2Pr[\text{no interruption at }\phi_1] \nonumber\\
    & +t_3y_3Pr[\text{no interruption at }\phi_1]Pr[\text{no interruption at }\phi_2] \nonumber\\
    & +t_4y_4Pr[\text{no interruption at }\phi_1]Pr[\text{no interruption at }\phi_2]Pr[\text{no interruption at }\phi_3] \nonumber\\
    & +\dots \nonumber\\
    =\  & t_1y_1+\sum_{i=2}^L\left(t_iy_i\prod_{k=1}^{i-1}Pr[\text{no interruption at }\phi_k]\right). \label{eq:f_rudimentaire}
\end{flalign}

An index~$k$ denotes a fidelity index smaller than or equal to a given~$i\in I$. The probability that no interruption happens at a sub-evaluation at~$\phi_k$ is the probability that all constraints assigned to~$\phi_k$ are satisfied after the sub-evaluation. Under the hypothesis that these probabilities are independent, the product of each~$(1-p_{kj})$ where~$B_{kj}=1$ for each~$j\in J$ is computed. Therefore,
\begin{equation}
    P_k(B):=Pr[\text{no interruption at }\phi_k]=\prod_{j\in J}(1-p_{kj}B_{kj}). \label{eq:Pno_int}
\end{equation}

The substitution of~\eqref{eq:Pno_int} in~\eqref{eq:f_rudimentaire} results in~\eqref{eq:modele_obj}. Problem~$\mathcal{Q}$ is defined as follows \\
\begin{minipage}{0.02\linewidth}
$\mathcal{Q}$
\end{minipage}
\begin{minipage}{0.9\linewidth}
\setlength{\belowdisplayskip}{0pt} \setlength{\belowdisplayshortskip}{0pt}
\begin{equation}
     \min_{B\in\mathbb{B}^{L\times m},y\in\mathbb{R}^L}\quad f(B)=t_1y_1+\sum_{i=2}^L t_iy_i\prod_{k=1}^{i-1}P_k(B) \label{eq:modele_obj}
\end{equation}
\begin{alignat}{2}
    \text{s.t.}\hspace{1cm} & \sum_{i\in I}B_{ij}=1\quad && \forall\,j\in J \label{eq:modele_part}\\
    & B_{ij}-\varepsilon\leq r_{ij}\quad && \forall (i,j)\in I\times J \label{eq:modele_tol}\\
    & B_{ij}\leq y_i\leq1\quad && \forall (i,j)\in I\times J \label{eq:modele_act1}\\
    & y_{i}\leq\sum_{j\in J}B_{ij}\quad && \forall\,i\in I. \label{eq:modele_act2}
\end{alignat}
\end{minipage} \\

Equation~\eqref{eq:modele_part} ensures that every blackbox constraint is assigned to exactly one fidelity, Equation~\eqref{eq:modele_tol} enforces the~$\varepsilon$ upper bound, and Equations~\eqref{eq:modele_act1} and~\eqref{eq:modele_act2} ensure the~$y_i$ variables respect their definition~\eqref{eq:def_y}.

Two assumptions are presumed by model~$\mathcal{Q}$. 
Although they are not required to apply the algorithm, 
    they are necessary to derive the theoretical results of Section~\ref{sec:reductions}.
\begin{assumption}
    The conditional sub-evaluation launched by Algorithm~\ref{algo:sous_eval_inter}'s second to last line has no impact on the optimal biadjacency matrix.
\end{assumption}
Since the moments when the extra sub-evaluations happen are unpredictable, model~$\mathcal{Q}$ disregards them.

\begin{assumption}
\label{eq:hyp_Q}
If a constraint is satisfied at its assigned fidelity, then it is also satisfied at any greater fidelity:
\begin{equation*}
        B_{ij}=1\text{ and }c_j(x,\phi_i)\leq0
        \quad \implies \quad p_{aj}=0\quad\forall\,a>i,\,\forall\,x\in X. 
    \end{equation*}
\end{assumption}
Define~$I_s=\{i\in I|e_i^\top B\neq0\}$ where $e_i$ is the $i$-th column of the $n\times n$ identity matrix. 
Without Assumption~\ref{eq:hyp_Q}, for each constraint~$c_j\leq0$, if~$B_{ij}=1$, then for each~$i_s\in I_s,i_s>i$, probability~$p_{i_sj}$ needs to be redefined as~$p_{i_sj}=Pr[c_j(x,\phi_{i_s})>0\,\lvert\,p_{i_s'j}\,\forall\,i_s'\in I_s,i\leq i'_s<i_s]$.
In practice, this assumption is sometimes violated, although in a proportion of cases of at most~$\varepsilon$ by definition of~$\varepsilon$.
Introducing Assumption~\ref{eq:hyp_Q} and selecting a small~$\varepsilon$ value is preferable not only because it simplifies the problem, but also because a small~$\varepsilon$ reduces the likelihood of misidentifying the feasibility of the constraints.

\subsection{Solving the model}
\label{sec:reductions}

Notice that Problem~$\mathcal{Q}$ is mixed-integer with a polynomial objective function. As today's solvers require a lot of solving time and have no guarantee of optimality for such problems, an alternative solving method is proposed. 
Indeed, the set of possible solutions may be reduced to the point where a simple exhaustive search is sufficient.

\subsubsection{Reduction by introducing non-differentiability}
\label{sec:red1}

The first reduction consists of simplifying the model by introducing non-differentiable elements, which are allowed since only an exhaustive search will be conducted.
First, the~$y_i$ variables are determined by the biadjacency matrix as:
\begin{equation*}
    y_{i}(B)=\left\{
    \begin{array}{ll}
        0\qquad & \text{ if }\sum\limits_{j\in J}B_{ij}=0 \\
        1 & \text{ otherwise}
    \end{array}
    \right.\quad\forall i\in I.
\end{equation*}

Hence, Equations~\eqref{eq:modele_act1} and~\eqref{eq:modele_act2} are removed, and a~$B$ matrix alone constitutes a solution to the model.
The function~$i:J\rightarrow I$ is introduced, which returns the index of the smallest fidelity in~$\Phi$ where constraint~$c_j\leq0$ can be assigned without violating Equation~\eqref{eq:modele_tol}:
\begin{equation*}
    i(j)=\text{min}\{i\in I:r_{ij}\geq1-\varepsilon\}.
\end{equation*}
Equation~\eqref{eq:modele_tol} is then equivalent to imposing that a constraint~$c_j\leq0$ can not be assigned to a fidelity~$\phi_i$ when~$i<i(j)$.
These new definitions allow for the definition of Problem~$\mathcal{Q}_1$, which is equivalent to~$\mathcal{Q}$.

\begin{minipage}{0.02\linewidth}
$\mathcal{Q}_1$
\end{minipage}
\begin{minipage}{0.9\linewidth}
\setlength{\belowdisplayskip}{0pt} \setlength{\belowdisplayshortskip}{0pt}
\begin{equation*}
    \min_{B\in\mathbb{B}^{L\times m}}\quad t_1y_1(B)+\sum_{i=2}^L\left(t_iy_i(B)\prod_{k=1}^{i-1}P_k(B)\right)\hspace{2.5cm}
\end{equation*}
\begin{alignat}{2}
    \text{s.t.}\hspace{1cm} & \sum_{i\in I}B_{ij}=1\quad && \forall\,j\in J \label{eq:modele_part2}\\
    & B_{ij}=0\quad && \forall\,(i,j)\in I \times J\mbox{ such that } i<i(j). \label{eq:modele_tol2}
\end{alignat}
\end{minipage}

\subsubsection{Reduction by filtering fidelities}
\label{sec:red2}

The second reduction aims to decrease the number of rows of the biadjacency matrix~$B$. It follows from Theorem~\ref{theo:reduction2}, which indicates that there exists an optimal solution for~$\mathcal{Q}_1$ where each constraint is assigned to a fidelity~$\phi_i$ such that
\begin{equation}
    i\in I_F:=\bigcup_{j\in J}i(j)\subseteq I, \label{eq:hatI}
\end{equation}
with~$I_F$ the filtered set of fidelity indexes. Each fidelity~$\phi_i$ where~$i\notin I_F$ can therefore be removed without excluding an optimal solution for~$\mathcal{Q}_1$, by replacing~$I$ with~$I_F$. Note that~$i(j)$ is a function that expresses a direct link between a fidelity index~$i$ and a constraint index~$j$. Conversely, the set~$I_F$ removes that link. Theorem~\ref{theo:reduction2} states that an optimal solution exists where all constraints are assigned to fidelities~$\phi_i$ where~$i=i(j)$ for some~$j\in J$, no matter what this~$j$ is. The theorem holds under the following two assumptions. 
In practice, they need not to be verified to apply the algorithm.
\begin{assumption}
    \label{eq:hyp}
    For any index~$j\in J$, the probabilities~$p_{ij}$ are monotonic decreasing with~$i$, for each~$i$ where constraint~$c_j\leq0$ can be assigned to fidelity~$\phi_i$ without violating Equation~\eqref{eq:modele_tol2}:
    \begin{equation*}
        p_{aj}\geq p_{bj}\quad\forall j\in J,\,\forall a,b\in I\text{ where }i(j)\leq a<b.
    \end{equation*}
\end{assumption}
This assumption simplifies the proof of Lemma~\ref{lem:optimality}. 
A more realistic and technical assumption in the context of blackbox optimization is discussed in~\cite{MScXL}, where it is shown that even when Assumption~\ref{eq:hyp} is violated, as~$\varepsilon$ becomes smaller, the likelihood of Theorem~\ref{theo:reduction2} being valid increases.

\begin{assumption}
    \label{eq:hyp_t}
    For any fixed value~$x\in X$,
    the function~$t(x,\phi)$ is monotonic increasing with respect to~$\phi$.
\end{assumption}
This implies that for a given set of LH samples,
\begin{equation*}
    t_{a}\leq t_{b}\quad\forall\,a,b\in I\text{ where }a<b.
\end{equation*}

\begin{lemma}
    \label{lem:feasibility}
    Let~$B$ be a feasible solution for~$\mathcal{Q}_1$.
    If there exists a fidelity index~$i'\in I\backslash I_F$ to which at least one blackbox constraint is assigned, then the matrix~$B'$ where
    \begin{equation}
        B'_{ij}=\left\{
        \begin{array}{ll}
            0\qquad & \text{ if $i=i'$} \\
            1 & \text{ if $i=i'-1$ and }B_{i'j}=1 \\
            B_{ij} & \text{ otherwise}
        \end{array}
        \right.\qquad\forall (i,j)\in I\times J \label{eq:Ba}
    \end{equation}
    is feasible for~$\mathcal{Q}_1$. \\
    \begin{proof}
        Let~$B$ be a feasible solution for~$\mathcal{Q}_1$ and $i'\in I\backslash I_F$ be a fidelity index to which at least one blackbox constraint is assigned. 
        Equation~\eqref{eq:modele_part2} is satisfied  by~$B'$, because it is satisfied by~$B$ and~$B'_{i'-1\,j}+B'_{i'j}=1=B_{i'-1\,j}+B_{i'j}$ for each~$j\in J$.
        Concerning Equation~\eqref{eq:modele_tol2}, on the one hand,
        if $(i,j) \in I \times J$ is such that $i<i(j)$, 
        then $B_{ij}=0$ since $B$ is feasible. 
        On the other hand, if $j'\in J$ is such that $B_{i'j'}=1$ then  $i'>i(j')$,
        which is equivalent to $i'-1\geq i(j')$. 
        The last condition in Equation~\eqref{eq:Ba} ensures that $B'_{ij} = B_{ij} = 0$ for each pair $(i,j)\in I\times J$ such that $i<i(j)$, 
        which implies that the matrix $B'$ is feasible for $\mathcal{Q}_1$.
    \end{proof} 
\end{lemma}

\begin{lemma}
    \label{lem:optimality}
    Let~$B$ be a feasible solution for~$\mathcal{Q}_1$. Under Assumption~\ref{eq:hyp}, if there exists a fidelity index~$i'\in I\backslash I_F$ to which at least one blackbox constraint is assigned, then the matrix~$B'$ given by~\eqref{eq:Ba} satisfies~$f(B')\leq f(B)$. \\
    \begin{proof}
        Let~$B$ be a feasible solution for~$\mathcal{Q}_1$ and $i'\in I\backslash I_F$ be a fidelity index to which at least one blackbox constraint is assigned. The objective function~\eqref{eq:modele_obj} maybe be divided into four terms, which correspond to the fidelity indexes smaller than~$i'-1$, equal to~$i'-1$, equal to~$i'$ and greater than~$i'$ of the sum
        $$f(B)=T_{<i'-1}(B)+T_{i'-1}(B)+T_{i'}(B)+T_{>i'}(B) $$
        where
        \begin{alignat*}{2}
            & T_{<i'-1}(B)\ = \ t_1y_1(B)+\sum_{i=2}^{i'-2}t_iy_i(B)\prod_{k=1}^{i-1}P_k(B), & \quad
            & T_{i'-1}(B)\ = \ t_{i'-1}y_{i'-1}(B)\prod_{k=1}^{i'-2}P_k(B),\\
            & T_{i'}(B)\ = \ t_{i'}y_{i'}(B)\prod_{k=1}^{i'-1}P_k(B), &
            & T_{>i'}(B)\ = \ \sum_{i=i'+1}^{L}t_iy_i(B)\prod_{k=1}^{i-1}P_k(B).
        \end{alignat*}
        Equation~\eqref{eq:Ba} ensures that the first term satisfies~$T_{<i'-1}(B)=T_{<i'-1}(B')$. Assumption~\ref{eq:hyp} ensures that the last term satisfies~$T_{>i'}(B)\geq T_{>i'}(B')$, as the only change from~$T_{>i'}(B)$ to~$T_{>i'}(B')$ is that all~$p_{i'j}$ become~$p_{i'-1\,j}$. 
        For the two central terms, two cases are considered. 
        First, if~$y_{i'-1}(B)=0$, then~$T_{i'-1}(B)=T_{i'}(B')=0$ as~$y_{i'}(B')=0$. 
        Furthermore, Assumption~\ref{eq:hyp_t} ensures that
        \begin{equation*}
            T_{i'}(B) 
                \ = \ t_{i'}\prod_{k=1}^{i'-1}P_k(B)
                \ = \ t_{i'}\prod_{k=1}^{i'-2}P_k(B')
                \ \geq \ t_{i'-1}\prod_{k=1}^{i'-2}P_k(B')
                \ = \ T_{i'-1}(B').
        \end{equation*}
        
        Second, if~$y_{i'-1}(B)=1$, then~$T_{i'-1}(B)=T_{i'-1}(B')$ as the sub-evaluations at~$\phi_{i'-1}$ and at lower fidelities are unchanged, and~$T_{i'}(B)>T_{i'}(B')=0$ as~$y_{i'}(B')=0$. In both cases, the sum of the four terms of~$f(B)$ is greater than or equal to those of~$f(B')$. Consequently,
        $f(B)\geq f(B')$.
        \end{proof} 
\end{lemma}

\begin{theorem}
    \label{theo:reduction2}
    There exists an optimal solution for~$\mathcal{Q}_1$ in which every blackbox constraint is assigned to a fidelity~$\phi_i$ where~$i\in I_F$. \\
    \begin{proof}
        The argmin of~$\mathcal{Q}_1$ can possibly contain only solutions where every blackbox constraint is assigned to a fidelity~$\phi_i$ where~$i\in I_F$. The theorem is then trivial. Otherwise, there exists a solution $B^*_0$ in the argmin of~$\mathcal{Q}_1$ such that the set
        \begin{gather*}
            I'(B^*_0):=\{i\in I\backslash I_F:y_i(B^*_0)=1\}
        \end{gather*}
        \noindent is not empty. 
        For a given~$i'\in I'(B^*_0)$, it is possible to find another solution,~$B^*_1$, where every constraint assigned to~$\phi_{i'}$ is rather assigned to~$\phi_{i'-1}$, and where every other assignment is unchanged. Lemma~\ref{lem:feasibility} indicates that~$B^*_1$ is feasible, and Lemma~\ref{lem:optimality} indicates that it yields an equal or better objective function value than~$B^*_0$. Therefore,~$B^*_1$ is also part of the~$\argmin$. Then, the set~$I'(B^*_1)$ can then be calculated. As long as~$I'$ is non-empty, this process is repeated to find~$B^*_2$,~$B^*_3$ and so on. The maximum number of such iterations is~$\max\{L-i(j):j\in J\}$, implying this process always terminates. When it does in~$K$ iterations,~$I'(B^*_K)=\varnothing$, and~$B^*_K$ is an optimal solution for~$\mathcal{Q}_1$ in which every blackbox constraint is assigned to a fidelity belonging to~$I_F$.
    \end{proof}
\end{theorem}

\subsubsection{Reduction by filtering constraints}
\label{sec:red3}

The third reduction simply consists of filtering out all a priori constraints, as they have no impact on the optimal biadjacency matrix. This reduces the number of columns of the biadjacency matrix~$B$. The set~$J$ is replaced for
\begin{equation}
    J_F:=\{j\in J:c_j\leq0\text{ is not an a priori constraint}\}, \label{eq:hatJ}
\end{equation}
the set of filtered constraint indexes.

\subsubsection{Exhaustive search}
\label{sec:red4}

First, consider~$\mathcal{Q}_1$, and replace~$I$ with~$ I_F$ and~$J$ with~$J_F$, therefore removing columns and rows of variables from matrix~$B$. 
This new problem is named~$\mathcal{Q}_2$. Equation~\eqref{eq:modele_tol} indicates that~$B_{ij}=0$ if~$r_{ij}<1-\varepsilon$.
Equation~\eqref{eq:modele_part} indicates that every constraint is assigned to exactly one fidelity.
The set of solutions for Problem~$\mathcal{Q}_2$ that satisfies these two equations is denoted~$\Omega_\mathcal{Q}$. 
For each solution in~$\Omega_\mathcal{Q}$, the objective function value is computed with~\eqref{eq:modele_obj}, and an optimal solution is found.
From this solution, an optimal biadjacency matrix~$B^*$ of size~$L\times m$ can be created by reintroducing columns~$J\backslash J_F$ and rows~$I\backslash I_F$, and by giving the value of~0 to these new elements.
Computing such an optimal assignment with the estimations~$\hat{p}_{ij}$,~$\hat{r}_{ij}$ and~$\hat{t}_{i}$ rather than the true values which are unavailable constitutes the second step of the {\em fidelity and interruption controlled optimization algorithm}.
Because points from~$\text{H}\cap\Omega^{ap}$ are used to compute these estimations, the a priori constraints are disregarded by model~$\mathcal{Q}$. This is not problematic because they have no impact on the optimal biadjacency matrix.

\subsection{Fidelity and interruption controlled optimization algorithm}
\label{sec:full_algo}

 The steps detailed in the previous sections are integrated in the {\em fidelity and interruption controlled optimization algorithm}, as shown in Algorithm~\ref{algo:hier_optimisation}.
Any blackbox optimization solver can be considered for the direct search step, including heuristic methods. This is possible since the interruption mechanism is not implemented in the solver,
  but as a wrapper around the blackbox.
From the solver's point of view, interruptions are part of the blackbox.
If an unconstrained optimization solver is chosen, we then recommend to use the extreme barrier function~$f_\Omega$.
The convergence guarantees of Algorithm~\ref{algo:hier_optimisation} are inherited from those of the solver.

\begin{algo}
\label{algo:hier_optimisation}
    \begin{center}
    \begin{table}[ht]
        \hypertarget{detailed_op_cstrs_hier}{}
        \centering \renewcommand{\arraystretch}{1.2}
        \begin{tabular}{L}
            \hline
            \textbf{Algorithm \thealgocounter: }\text{Fidelity and interruption controlled optimization algorithm} \\
            \hline
            \textbf{Input:} \\
            \left.\middle|
               \begin{array}{rl}
                    x^0 & \text{optimization starting point (optional)} \\
                    \eqref{eq-pb} & \text{problem containing $X \subseteq \R^n$, $f$ and $c$} \\
                    \Phi & \text{ordered set of fidelities ending with $1$} \\
                    \varepsilon & \text{upper bound on the probability that a constraint's feasibility is misidentified} \\
                    \rho, n_H & \text{Latin hypercube sizing factor and sample size} \\
                    \textsf{solver} & \text{direct search blackbox optimization solver} 
                \end{array}
            \right. \\
            \textbf{1. }\text{Constraint feasibility relative to fidelity model.} \\
            \left.\middle|
               \begin{array}{l}
                    \textbf{If }x^0\text{ is provided, find the $(\ell^\text{cen},u^\text{cen})$ bounds from $X$, $\rho$ and $x^0$ using~\eqref{eq:licen} and~\eqref{eq:uicen}}. \\
                    \textbf{Else}\text{, set }(\ell^{cen},u^{cen})=(\ell,u)\text{ and at end of step 1, set $x^0$ as the best point in $H$.} \\
                    \text{Randomly determine $H$, the $n_H$ LH points bounded by $\ell^{cen}$ and $u^{cen}$.} \\
                    \text{Evaluate each point in $H$ at each fidelity in $\Phi$ by parallelizing as much as possible.} \\
                    \text{Calculate all $\hat{r}_{ij}$, $\hat{p}_{ij}$ and $\hat{t}_i$
                    estimations using~\eqref{eq:hatr},~\eqref{eq:hatp} and~\eqref{eq:hatt}, respectively.}
                \end{array}
            \right. \\
            \textbf{2. }\text{Optimal biadjacency matrix computation.} \\
            \left.\middle|
               \begin{array}{l}
                    \text{Find $J_F$, and find $ I_F$ with $\varepsilon$
                    using~\eqref{eq:hatJ} and~\eqref{eq:hatI}, respectively}. \\
                    \text{Solve problem $\mathcal{Q}_2$ with an exhaustive search on }\Omega_\mathcal{Q}. \\
                    \text{Create a matrix $B\in\mathbb{B}^{L\times m}$ from an optimal solution for $\mathcal{Q}_2$ using $\hat{r}_{ij}$, $\hat{p}_{ij}$ and $\hat{t}_i$}. \\
                \end{array}
            \right. \\
            \textbf{3. }\text{Direct search.} \\
            \left.\middle|
               \begin{array}{l}
                    \text{Initialise $f^*=\infty$.} \\
                    \text{Launch \textsf{solver}, providing Algorithm~\ref{algo:sous_eval_inter} with matrix $B$ as the blackbox evaluation} \\
                    \quad\text{function, and $x^0$ as the starting point if needed. Update $f^*$ after each evaluation.} \\
                \end{array}
            \right. \\
            \textbf{Return }\text{the \textsf{solver} output} \\
            \hline
        \end{tabular}
    \end{table}
    \end{center}
\end{algo}

\clearpage
\newpage
\section{Computational results}
\label{sec:results}

The interruption mechanism described above is tested on four instances of the 
\solar{}\footnote{Available at
\href{https://github.com/bbopt/solar}{\url{https://github.com/bbopt/solar}} (version~1.0)}
family of blackbox problems~\cite{solar_paper}.
To the best of the authors' knowledge, \solar{} is the only benchmark blackbox problem collection where constraints can depend on fidelity.
The properties of the four \solar{} instances used herein are presented in Table~\ref{tab:solar_instances}.

\begin{table}[ht!]
    \renewcommand{\tabcolsep}{4.5pt}
        \begin{tabular}{|C|CCCC|}
        \hline
        \text{Instance} & n & m & \text{nb. a priori constraints} & \text{nb. multi-fidelity constraints} \\
        \hline\hline
        \text{\solar{2}} & 14 & 13 & 5 & 4 \\
        \hline
        \text{\solar{3}} & 20 & 13 & 5 & 5 \\
        \hline
        \text{\solar{4}} & 29 & 16 & 7 & 6 \\
        \hline
        \text{\solar{7}} & 7 & 6 & 2 & 2 \\
        \hline
        \end{tabular}
    \caption{Characteristics of the four studied blackbox problem instances from the \solar{} family of problems.}
    \label{tab:solar_instances}
\end{table}

The \nomad V4~\cite{nomad4paper} software package is chosen to quantify the effect of the interruption mechanism from Algorithm~\ref{algo:hier_optimisation}.
\nomad is an implementation of the \mads algorithm~\cite{AuDe2006, AuDe09a}
    conceived for constrained blackbox optimization problems.
\nomad is used because that it has shown to be successful on real engineering and industrial problems~\cite{AlAuGhKoLed2020}, it is freely available and can easily be modified to include the interruption mechanisms proposed here.
Comparisons with other blackbox solvers, and tuning the \nomad parameters are beyond the scope of the present paper.

With \nomad, constraints can be managed using the EB or the PB strategies.
Literature suggests that  PB is preferable~\cite{AuDe09a}.
However, due to Algorithm~\ref{algo:sous_eval_inter} not returning the true outputs when the evaluation of a point is interrupted, the EB might be more adapted, as it rejects these points.
Conversely, the PB uses output values to compute new incumbent points.
Therefore, two implementations and a base case are tested for comparison:
\begin{itemize}
    \item Inter PB: Algorithm~\ref{algo:hier_optimisation} with \nomad and the PB;
    \item Inter EB: Algorithm~\ref{algo:hier_optimisation} with \nomad and the EB;
    \item Base case: \nomad with default parameters and blackbox fidelity fixed at~1.
\end{itemize}

The method introduced in this paper is motivated by computationally expensive problems, such as Hydro-Qu\'ebec's asset management blackbox optimization problem from PRIAD, where the solver and the Algorithm~\ref{algo:sous_eval_inter} computation times are insignificant in comparison. To replicate such problems with \solar{} instances that  are computationally less demanding but possess desirable characteristics, only blackbox computation times are considered in the following profiles. For both implementations of the algorithm, the empirically determined values of $n_H=10^4$, $\varepsilon=0.05$ and
\begin{equation*}
    \Phi=\{0,10^{-10},2^{-10},2^{-9},2^{-8},2^{-7},2^{-6},2^{-5},2^{-4},0.1,0.2,0.3,0.4,0.5,0.6,0.7,0.8,0.9,1\}
\end{equation*}
are chosen. All numerical experiments are performed using an Intel Xeon Gold 6150 CPU @ 2.70GHz processors.

\subsection{Without a starting point for the optimization}
\label{sec:no_x0}

\begin{figure}[htb]
    \centering
    \begin{subfigure}[t]{0.49\textwidth}
    \includegraphics[width=\textwidth]{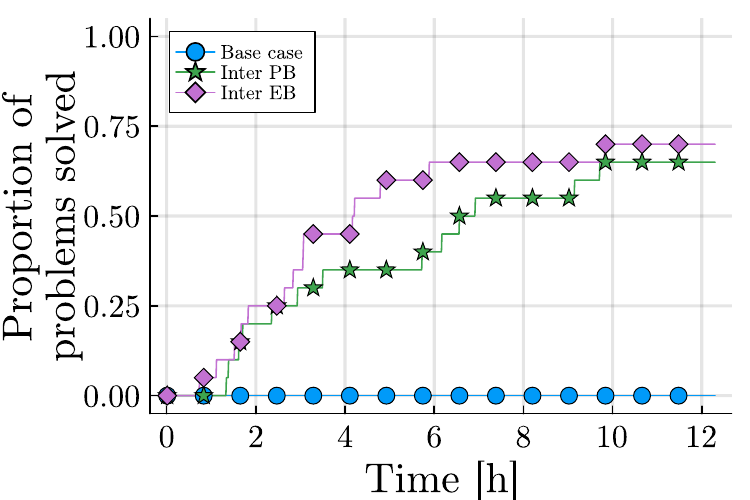}
        \caption{$\tau=0.05$}
        \label{fig:data_prof_hl_part_fea_start_s2_0.05.pdf}
    \end{subfigure}
    \hfill
    \begin{subfigure}[t]{0.49\textwidth}
        \includegraphics[width=\textwidth]{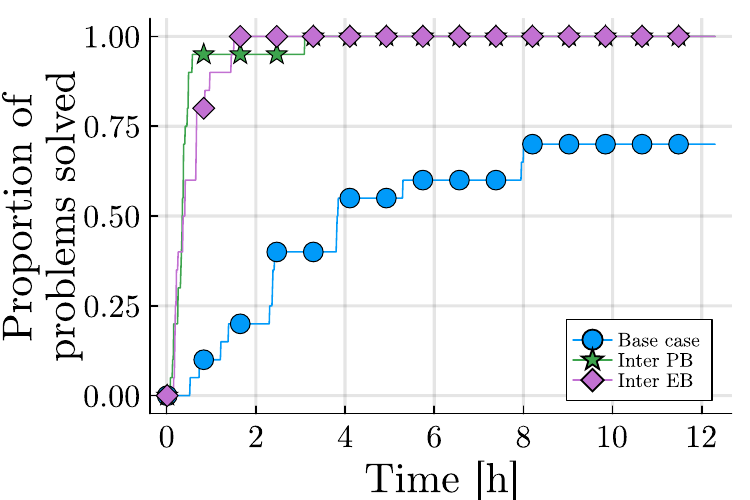}
        \caption{$\tau=0.5$}
        \label{fig:data_prof_hl_part_fea_start_s2_0.5.pdf}
    \end{subfigure}
    \caption[Profils de donn\'ees de \solar{2} à partir de points r\'ealisables.]{\solar{2} data profiles from 20 runs with no $x^0$.}
    \label{fig:data_prof_hl_part_fea_start_s2}
\end{figure}

This section presents optimizations where no starting point is provided by the user. Consequently, the LH bounds are those of~$X$. The base case also performs a LH with the same parameters to find a starting point. Optimizations on three constrained multi-fidelity instances of the \solar family of blackboxes: \solar{2}, \solar{3} and \solar{4} are conducted. By varying the feasible starting point, 20 optimization runs are executed for each tested instance. The results are illustrated in data profiles, with the initial two profiles relating to \solar{2}, and shown in Figure~\ref{fig:data_prof_hl_part_fea_start_s2}. As the LH times are identical, only the optimization times are shown.

The Algorithm~\ref{algo:hier_optimisation}'s success, whether paired with the EB or the PB, is mainly attributed to two factors.
Firstly, \solar{2} contains a frequently violated constraint with a high estimated probability of representativity at low fidelities, which allows for numerous quick interruptions on infeasible points.
Secondly, the calculated~$B$ matrix is accurate because the constraint feasibility relative to fidelity model is accurate for most points encountered during optimization.
The opposite is true for \solar{3} and \solar{4}.
This happens because the constraints' behaviour is non-homogeneous throughout~$X$, that is, the~$r_{ij}$,~$p_{ij}$ and~$t_i$ values vary greatly when they are calculated from different sub-spaces of~$X$.
As a result, when using Algorithm~\ref{algo:hier_optimisation} with those two instances, every evaluated point is systematically estimated to be infeasible. No figures are shown as there are no curves for the~\ref{algo:hier_optimisation} implementations. Conversely, the base case finds numerous feasible solutions.

\subsection{With a starting point for the optimization}
\label{sec:yes_x0}
\begin{figure}[htb]
    \centering
    \begin{subfigure}[t]{0.49\textwidth}
        \includegraphics[width=\textwidth]{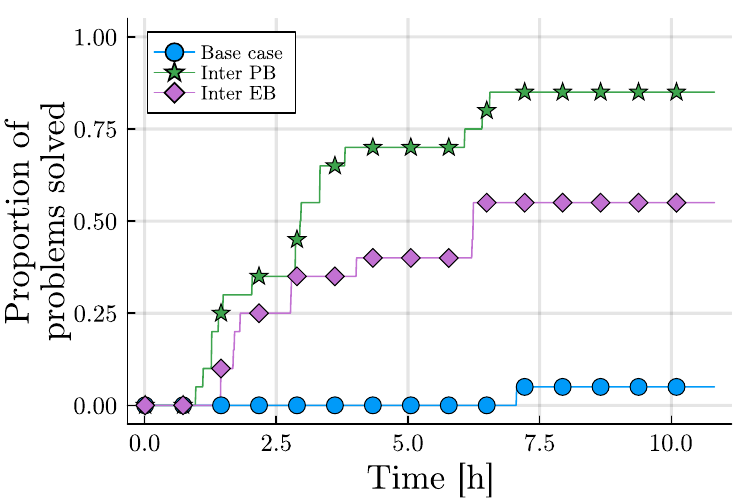}
        \caption{$\tau=0.1$}
        \label{fig:data_prof_cen_hl_fea_start_s2_0.1}
    \end{subfigure}
    \hfill
    \begin{subfigure}[t]{0.49\textwidth}
    \includegraphics[width=\textwidth]{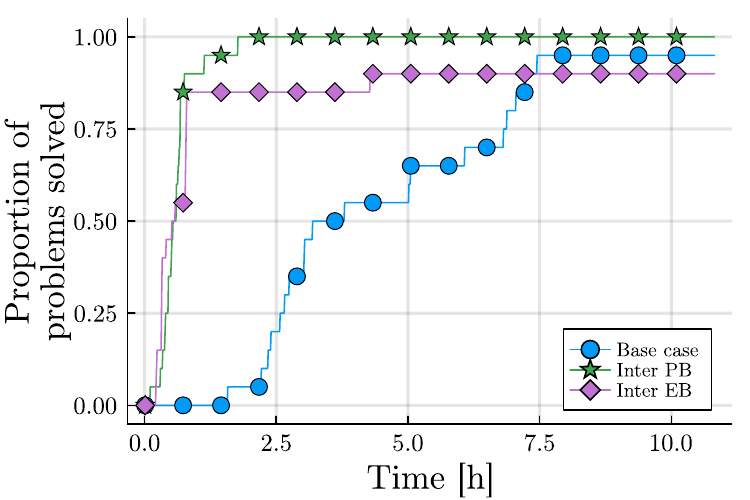}
        \caption{$\tau=0.5$}
        \label{fig:data_prof_cen_hl_fea_start_s2_0.5}
    \end{subfigure}
    \caption[Profils de donn\'ees de \solar{2} avec un hypercube latin centr\'e.]{\solar{2} data profiles with a starting point. Curves from implementations of Algorithm~\ref{algo:hier_optimisation} start at 642.22 seconds to account for LH sampling time.}
    \label{fig:data_prof_cen_hl_fea_start_s2}
\end{figure}

This section presents optimizations where a known feasible starting point is provided by the user.
The optimization runs are conducted on each constrained multi-fidelity instances of the \solar family of blackboxes: \solar{2}, \solar{3}, \solar{4} and \solar{7}, in which $\rho$ is set to~$\frac{1}{4}$,~$\frac{1}{10}$,~$\frac{1}{20}$ and~$\frac{1}{4}$ respectively.
Those values are based on preliminary results.
In this scenario, the base case does not perform any LH sampling, which grants it a time advantage.
However, this advantage is inconsequential due to the extensive parallelization of samples using Hydro-Qu\'ebec's facilities.
By varying the \nomad seed for random polls, 20 optimization runs are executed for each tested instance, and the results are illustrated in data profiles.
Two data profiles for \solar{2} are shown in Figure~\ref{fig:data_prof_cen_hl_fea_start_s2}.

With~$\tau=0.5$, the base case solves a greater number of problems compared to the proposed algorithm with the EB, as the implementation is highly inefficient for one of the 20 optimization runs.
In general, both implementations of Algorithm~\ref{algo:hier_optimisation} are preferable.

\begin{figure}[!htb]
    \centering
    \begin{subfigure}[t]{0.49\textwidth}
        \includegraphics[width=\textwidth]{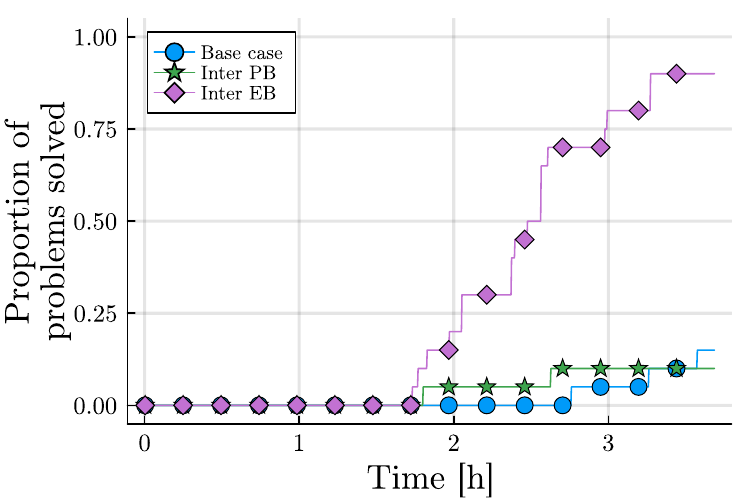}
        \caption{$\tau=0.05$}
        \label{fig:data_prof_cen_hl_fea_start_s3_0.05}
    \end{subfigure}
    \hfill
    \begin{subfigure}[t]{0.49\textwidth}
    \includegraphics[width=\textwidth]{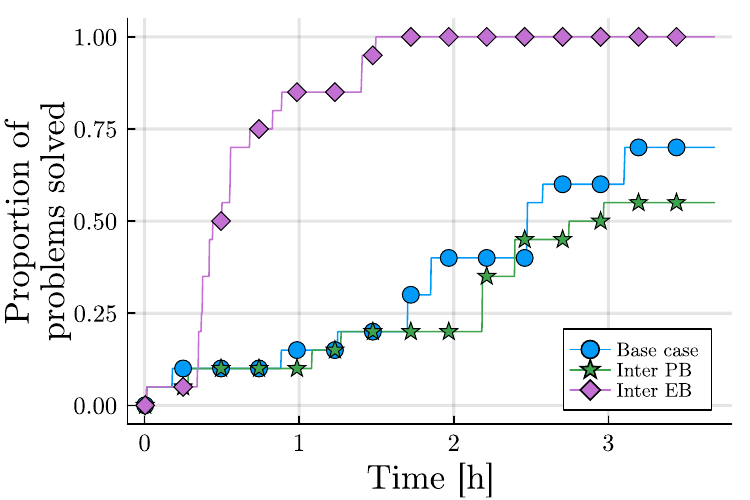}
        \caption{$\tau=0.3$}
        \label{fig:data_prof_cen_hl_fea_start_s3_0.3}
    \end{subfigure}
    \caption[Profils de donn\'ees de \solar{3} avec un hypercube latin centr\'e.]{\solar{3} data profiles from 20 runs with a given $x^0$. Curves from implementations of Algorithm~\ref{algo:hier_optimisation} start at 45.64 seconds to account for LH sampling time.}
    \label{fig:data_prof_cen_hl_fea_start_s3}
\end{figure}

\begin{figure}[htb]
    \centering
    \begin{subfigure}[t]{0.49\textwidth}
        \includegraphics[width=\textwidth]{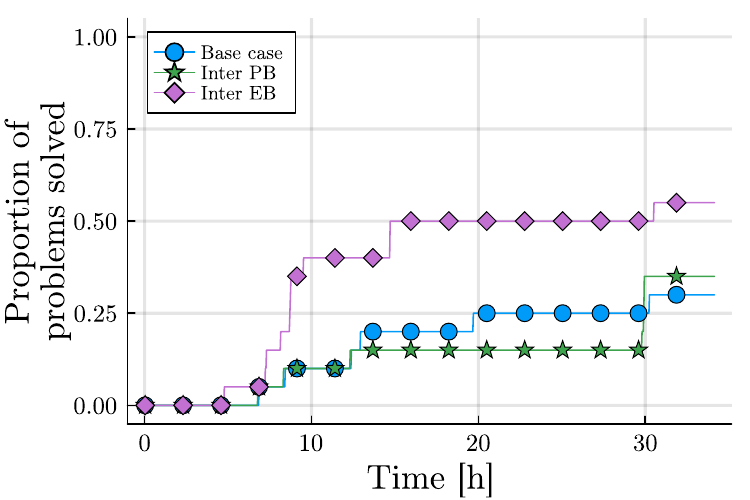}
        \caption{$\tau=0.01$}
        \label{fig:data_prof_cen_hl_fea_start_s4_0.01}
    \end{subfigure}
    \hfill
    \begin{subfigure}[t]{0.49\textwidth}
    \includegraphics[width=\textwidth]{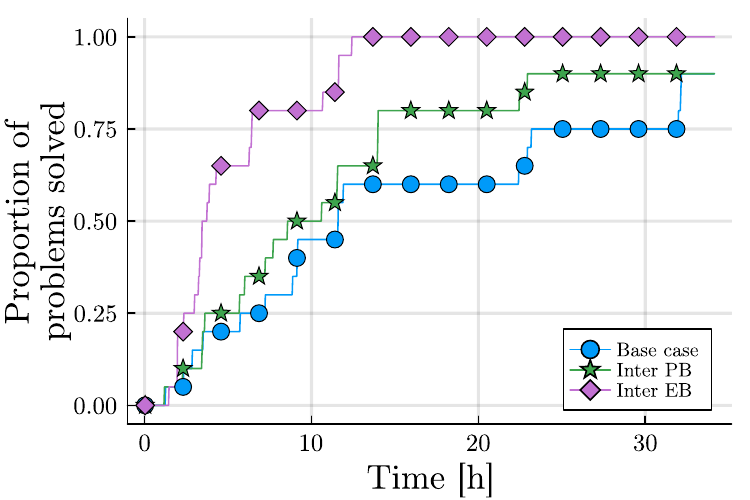}
        \caption{$\tau=0.1$}
        \label{fig:data_prof_cen_hl_fea_start_s4_0.1}
    \end{subfigure}
    \caption[Profils de donn\'ees de \solar{4} avec un hypercube latin centr\'e.]{\solar{4} data profiles from 20 runs with a given $x^0$. Curves from implementations of Algorithm~\ref{algo:hier_optimisation} start at 94.75 seconds to account for LH sampling time.}
    \label{fig:data_prof_cen_hl_fea_start_s4}
\end{figure}

Data profiles for \solar{3} and \solar{4} are shown in Figure~\ref{fig:data_prof_cen_hl_fea_start_s3} and Figure~\ref{fig:data_prof_cen_hl_fea_start_s4}, respectively.

For both instances, the base case yields results comparable to Algorithm~\ref{algo:hier_optimisation} paired with the PB.
On the other hand, when pairing the algorithm with the EB, it performs significantly better. This can be attributed to a higher scarcity of feasible points in \solar{3} and \solar{4} compared to \solar{2}.
This section also studies \solar{7}, a constrained multi-fidelity blackbox where infeasible points are much less common than for the other tested instances.
It serves as a test to assess how Algorithm~\ref{algo:hier_optimisation} performs when it has limited opportunities to interrupt evaluations and save time in contrast to the base case.
Additionally, the objective function value is affected by multi-fidelity for this instance.
Thus, the condition~$y_{L}=1$ is imposed.
Results show that the optimal biadjacency matrix computed for this problem assigns all constraints to~$\phi_L=1$.
This suggests that the method has discerned the absence of meaningful opportunities for interruptions and that emulating the base case is the optimal approach.
No figure is shown; the base case and the PB implementation exhibit identical data profiles, except for the fact that the base case consistently precedes by 181.07 seconds due to its absence of LH sampling prior to any optimization.


\section{Discussion}
\label{sec:discu}

This work introduces a novel approach to computationally expensive multi-fidelity black\-box optimization problems by leveraging low-fidelity assessments of constraints violation to interrupt evaluations. These assessments are determined by a biadjacency matrix, which balances computation cost and the probability of constraint violation. Our computational results demonstrate that, under specific conditions, pairing the \nomad solver with the {\em fidelity and interruption controlled optimization algorithm} yields significantly superior solutions compared to \nomad alone. Here is a summary of favorable conditions:
\begin{itemize}
    \setlength\itemsep{0.3em}
    \item scarce feasible points;
    \item accurate constraint violation assessments at lower fidelities;
    \item homogeneity in constraint behaviour relative to fidelity within the LH bounds defined by the sizing factor~$\rho$.
\end{itemize}

 This final condition is necessary for an accurate constraint feasibility relative to fidelity model, which results in a high quality biadjacency matrix.
When this condition is not fulfilled with~$\rho=1$, the existence of a known feasible solution prior to the optimization becomes vital; it enables the selection of a sizing factor~$\rho$ that increases the homogeneity in constraint behaviour.

 While the scarcity of feasible points typically presents challenges, the proposed method achieves greater computational cost reductions under such a condition.

When utilizing the \nomad solver, we observe that the preferred barrier choice depends on the blackbox. For problems with infrequent feasible points, the EB is more suitable, while the PB is mostly preferred when feasible points are more common.

Future work involves dynamically computing the biadjacency matrix and updating the constraint feasibility relative to fidelity model, improving this model by introducing machine learning methods such as~\cite{AuCaJa2020a,MeMuSi2023} as well as tools from trust-region methods and developing a multi-fidelity barrier method.
Additionally, the proposed algorithm will be applied to complex industrial problems such as PRIAD, and adapted to these problems when an exploitable structure from the blackbox is available.

\backmatter
\section*{Data availability and conflict of interest statement}

The solar problems are available at \href{https://github.com/bbopt/solar}{\url{https://github.com/bbopt/solar}}.
The NOMAD software package is available at \href{https://github.com/bbopt/nomad}{\url{https://github.com/bbopt/nomad}}.
This work has been funded by
the NSERC Discovery grants
of C.~Audet~(\#2020-04448)
and S.~Le~Digabel~(\#2018-05286),
by the NSERC/Mitacs Alliance grant~(\#571311-21) in collaboration with Hydro-Qu\'ebec,
and by X.~Lebeuf's NSERC Graduate Scholarships~–~Master’s grant.

\end{document}